\documentclass[PRL ]{revtex4-1}

\usepackage{graphicx}
\usepackage{epstopdf}
\usepackage{mathptmx}      
\usepackage{subfigure}
\usepackage{amsmath}
\usepackage{algorithm}
\usepackage{algorithmic}
\usepackage{multirow}
\usepackage{url}
\usepackage{natbib}

\begin{document}
\title{The asymptotic behaviors of self excitation information diffusion processes for a large number of individuals}
\author{Lifu Wang}
\email[]{Lifu\_Wang@bjtu.edu.cn}
\affiliation{Department of Electronic and Information Engineering, Beijing Jiaotong University, Beijing, China
Key Laboratory of Communication and Information Systems, Beijing Municipal Commission of Education Beijing Jiaotong University, Beijing, China}
\author{Bo Shen}
\email[]{bshen@bjtu.edu.cn}
\affiliation{Department of Electronic and Information Engineering, Beijing Jiaotong University, Beijing, China
Key Laboratory of Communication and Information Systems, Beijing Municipal Commission of Education Beijing Jiaotong University, Beijing, China}
\date{\today}
\begin{abstract}
The dynamics of opinion is a complex and interesting process, especially for the systems with large number individuals. It is usually hard to describe the evolutionary features of these systems. In this paper, we study the self excitation opinion model, which has been shown the superior performance in learning and predicting opinions. We study the asymptotic behaviors of this model for large number of individuals, and prove that the asymptotic behaviors of the model in which the interaction is a multivariate self excitation process with  exponential function weight, can be described by a Mckean-Vlasov type integro differential equation.  The coupling between this equation and the initial distribution captures the influence of self excitation process, which decribes the mutually- exicting and recurrent nature of individuals. Finally we show that the steady state distribution is a ``contraction" of the initial distribution in the linear interaction cases.

\end{abstract}
\pacs{}
\maketitle

\section{Introduction}
\label{intro}
Expressing opinions and then influencing others are the main forms of social behavior of people. The dynamics of public opinion has long been a major concern of social science and  computational social science. There has been an increasing interest on the perception and prediction of social opinion, such as quantitative investment firms, who measure investor sentiment and trade using social media\cite{Karppi2016Social}, and prediction of election results\citep{ceron2015using}.

Based on sentiment analysis with deep learning, there are many frameworks to perceive opinion on the social network. For the prediction of social opinion,
there are some noticeable models  proposed at the beginning of the 21st century, such as Sznajd model\citep{Katarzyna2000Opinion}, Hegselmann-Krause model \citep{Hegselmann2005Opinion, Hegselmann2002Opinion}, Deffuant-Weisbuch model \citep{Guillaume2000Mixing}, and so on, which had been established to understand the coordinated movements of opinion as a group. However, these work held some limitations: (i) Most of these model on opinion dynamics are theoretical and have not proved their effectiveness quantificationally. (ii) These models are represented by a cellular automata, which makes it difficult to analyze the global behavior of the dynamics.

More recently, some researchers began to study the models which can provide more accurate predictions\citep{Das2014Modeling,De2014Learning,De2015Learning}. The authors in \citep{De2015Learning} proposed a framework of opinion dynamics named SLANT, providing  accurate predictions of users' opinions. In the early models, the opinion and the states of the next moment only depend on the states of the last moment, thus have nothing to do with more previous ones. It means that opinion formation is seen as a Markov process. However, in the real world, this assumption may not be precise. The model in \citep{De2015Learning} indicates that self excitation opinion models with  exponential function weight have superior performance in learning and predicting opinions.

In many works \citep{Zhou2013Learning,De2015Learning,Wang2016A,wang2017variational,Wang2017}, self excitation models have shown efficient abilities to capture users' behavior pattern. However, these models still base on cellular automata, thus can not be employed to characterize the global dynamics of the opinion evolution analytically.

In this paper, we study the self excitation opinion dynamics on the large homogeneous network. We show that when the number of the individuals is very large, the opinion distribution of this system evolve according to a  nonlinear partial differential equation(PDE) of Mckean-Vlasov type. Due to the dependence on the initial value in the local rules, the individuals with different initial value will have different dynamic equations, thus the the mean field equation deriving process is different from the cases of the previous works\citep{M1996Asymptotic,Andreis2017McKean}. We show that this PDE can be decomposed into two part, where  the first part is the same as the Fokker-Planck model in \citep{Toscani2006Kinetic}. And the second part  can capture the mutually-exciting and recurrent nature of individual behaviors by an integral term related to the initial value. This term leads to different steady state behaviors from the Markov model such as \citep{Toscani2006Kinetic}. Generally, the stationary solutions  of Toscani type model will trend to a single-point support Dirac measure when noise parameters tend to 0, but the model with additional term makes the steady state to be a ``contraction" of the initial distribution.

The rest of this paper is organized as follows. Section 2 is for related work. In section 3, we introduce the opinion dynamics we study and the outline of our proof. Some examples of model and the analysis are shown in section 4, and section 5 is a summary. The detail proof is in the appendix.

\section{Related work}

\subsection{Information diffusion model}
The early theory of information diffusion in mathematical sociology was proposed by Sznajd in \citep{Katarzyna2000Opinion}. Sznajd's model is designed to explain the features of opinion dynamics, in which every individual is on a lattice and have a state 0 or 1 to express the opinion. Individuals update their opinion by the opinions of their neighbours. The numerical simulations of this model are also investigated by many researchers \citep{Bernardes2002Election,Dietrich2002Better,Stauffer2002Persistence}. On the complete graph, Slanina \citep{Slanina2003Analytical} shows that when the number of individuals is very large, the probability densities evolve according to the partial differential equations
\begin{equation}\label{a1}
\frac{\partial}{\partial \tau}P(m,\tau)=\frac{\partial^2}{\partial m^2}[(1-m^2)P(m,\tau)]
\end{equation}
And
\begin{equation}\label{a2}
\frac{\partial}{\partial \tau}P(m,\tau)=-\frac{\partial^2}{\partial m^2}[(1-m^2)mP(m,\tau)]
\end{equation}
where $P(m)$ is the probability distribution of opinion. Equation (\ref{a1}) describes the Ochrombel simplification of Sznajd model, and (\ref{a2}) does the original case. The authors showed the existence of phase transition in the original formulation and smooth behaviour without transition in the model modified by Ochrombel.

For the theory of continuous opinion dynamics, Deffuant-Weisbuch model \citep{Guillaume2000Mixing} and Hegselmann-Krause model\citep{Hegselmann2005Opinion,Hegselmann2002Opinion} deal with two different cases. Both these models are so called bounded confidence models, i.e. the individuals can only be influenced by some others whose opinions are close enough to theirs. The characteristics of bounded confidence lead to the formation of communities, which coincides with the phenomenon in the real world. A very beautiful analysis of the convergence of Deffuant-Weisbuch model is given in \citep{Zhang2012Convergence}, which reveals the community structure of the model. Ref.\citep{Bhattacharyya2013On} proves the convergence of Hegselmann-Krause model. On the other hand, numerical simulations on complex network \citep{Guo2009Bifurcation} show that ``the more heterogeneous the complex network is, the weaker the ability of polarization and consensus of complex network will be.''. Models with both continuous opinion and discrete states are also considered, such as SHIR model in \citep{Liu2016SHIR}, where the individuals have a four-states variable to represent susceptible, hesitated, infected and removed.

In \citep{Como2011SCALING}, the authors gived a general framework to prove that the scaling limits of many-body continuous opinion system can be described by  a measure value ordinary differential equation(ODE). In such models, every individual will change its opinion to $\omega X +(1-\omega)X'$, a convex combination with the neighbor with probability $k(X,X')$. When $k(X,X')=1_{[0,\epsilon]}|(|X-X'|)$, it becomes Deffuant-Weisbunch model. The key ideas of their work are the large deviation principle and some estimations based on the Lipschitz property, which are the same as the general method used to prove the chaotic propagation properties \citep{Sznitman1991Topics}.

Based on  Boltzmann equation with granular gas like interactions, Toscani\citep{Toscani2006Kinetic} introduced a collision model, in which the communication between individuals was considered as a collision with the following form,
\begin{equation}
w'=(1-\eta P(w,w_*))w+\eta P(w,w_*)w_* +\zeta D(w)
\end{equation}
where function D and P describe the local relevance of the compromise and diffusion for a given opinion. Toscani also considered the quasi invariant opinion limit for Boltzmann type equation. The main idea of that is to scale the interaction frequency, strength and the diffusion in the integral equation, then the equation reduces to a Fokker-Planck type equation. Then on the base of this model, considering the influence of the structure of social network, the authors in \citep{Albi2016Opinion} consider the kinetic opinion dynamics on the large scale networks evolving over time. The evolutions of both the network and the opinion were involved.

\subsection{Self excitation point process in social systems}
Point processes type models are generally used to analyze the impact of events on the system such as CSMA/CD\citep{Cho2010A}, financial data\citep{Chavez2012High} and the dynamics of book sales\citep{Desch2005Dynamics}. A very succinct and effective analysis for the time series of daily views for videos on YouTube is introduced in \citep{Crane2008Robust}. The authors showed that, most of the video viewing record data can be described statistically by a Poisson process, but still about 10\% data show a point process with instantaneous rate $\lambda(t)$
\begin{equation}
\lambda(t) =V(t)+\sum_{i,t_j\leq t} \mu_i \phi(t-t_i)
\end{equation}
where $\mu_i$ is the number of potential viewers influenced directly by person $i$ who views the video through social network, and $V(t)$ captures all spontaneous views that are not triggered by epidemic effects on the network.

For the self excitation social opinion models, the authors in \citep{De2015Learning} proposed a opinion model on a social network, in which every user has a latent opinion about the given topic and can post messages about the topic. The experiments on real data show that self excitation  model, with the intensity depending on the messages sent from a neighbor in the past and the weight function having the form like $e^{-\omega t}$, performs much more better than the Poisson process($\omega=0$) model.

In \citep{Wang2017}, the authors used point process to predicting user activities. They proposed a generic framework  for point process prediction problem and they used a mass transport equation to update  the transport rate and compute the conditional mass function.

\subsection{Mean field theory}

Mean field theroy(MFT) is introduced to study many-body problem by using a single averaged effect to approximate the effect of all individuals. The core of MFT is to estimate the error of the mean field approximation. If there is no interaction, all individuals are independent of each other, then the  law of large numbers works and gives the mean field approximation. In the general case, we still want to be able to use the law of large numbers. That is to say, if one picks a chaotic (i.i.d) initial distribution of particles, we hope that this distribution is still chaotic as time evolving, which is the so called propagation of chaos.

The theory of propagation of chaos in the cases of Wiener noise is summed  in \citep{Sznitman1991Topics}. if $N$ particles with initial "chaotic" distribution $u_0^{\otimes N}$ satisfy the  stochastic differential equation(SDE)
\begin{equation}
dX_i=dW_i+\frac{1}{N}\sum_j b(X_i,X_j)dt
\end{equation}
the Mckean-Vlasov mean field equation is
\begin{equation}
dX=dW+\int b(X,y)u(dy) dt
\end{equation}
where u is the law of $X$. The method to show this is to use the Lipschitz character of $b$ and Gronwall's lemma. Mckean-Vlasov equation for SDE with Possion jumps is in \citep{Andreis2017McKean}. The authors used Burkholder-Davis-Gundy inequality for martingales. The main step in \citep{Andreis2017McKean} is to use Doob-Meyer decomposition and  estimate the compensated Poisson process.

Using the large deviation principle \citep{Deuschel1989Large}, Arous and Guionnet\citep{Arous1995Large} studied the mean field simplification dynamics for Langevin spin glass. The framework of their proof makes use of Sanov's Theorem  and Varadhan's Lemma in the path space $C[0,T]$ to get the rate function for the interaction system. This method was generalized to the interacting stochastic processes in random media in \citep{Pra1996McKean}.

The main obstacles to use these methods to prove MFT are the unboundedness of coefficients in SDE. Because of Varadhan's lemma, if we use a stop time to get the localization of SDE, we can get the so called local large deviation principle. Dawsont and Gartner \citep{Dawsont1987Large} gived some compactness criterion to convert the "local" result into "global" one. They also show that a Lyapunov function for the system of weakly interacting diffusions will  let the compactness condition be satisfied. Puhalskii \citep{Puhalskii2001Large} set up the whole framework of local to global LDP and introduce the C-exponential tightness conditions.

\section{The  dynamics and the framework of the proof}
\subsection{Self excitation  interaction system }
The model proposed in \citep{De2015Learning} has the following form.\\
Opinion model:
\begin{equation}
 x_i =b_i +\sum_j \alpha_{i,j} \sum_{t_j\in H_j(t)} m_j k(t-t_j)   =b_i +\sum_j \alpha_{i,j} m_j k(t)*dN_j(t)
\end{equation}
where  $m_j$ is the message send by user j, $H_j(t)$ is the history of events up to time t for user j, "*" is the convolution operator, $k(t)$ is a triggering kernel and N(t) is a Possion process with intensity $\lambda$.\\
Messages intensity model:
For the intensity to send a message $\lambda(t)$, we have
\begin{equation}
\lambda(t)_i =\mu_i +\sum_j b_{i,j} g(t)*dN_j(t)
\end{equation}
And the values of message m come from a sentiment distribution related to the opinion x and $E[m_i|x_i]=x_i$.\\
Then it can be simplified if we only want to consider opinion part\citep{Wang2016A,wang2017variational} :
\begin{equation}
 x_i =b_i +\sum_j h(x_j,x_i) k(t)*dN_j(t)
\end{equation}

In the case that $h(x_i,x_j)=h(x_j)$ where $x_j$ is the opinion of the individual $j$, ref.\citep{De2015Learning,Wang2016A} showed the validity of the model on the real world data. It is also reasonable to consider the case that $h(x_j,x_i)=x_j-x_i $, so that the form of this model becomes the same as the case in DW model\citep{Guillaume2000Mixing} and Toscani model \citep{Toscani2006Kinetic}.

In \citep{De2015Learning,Wang2016A}, self excitation opinion model with exponential function type weight was shown to be very  effective in the opinion dynamics, therefore we set $k(t)=e^{-\omega t}$, where $\omega$ is a constant. Due to the fact that differential of the convolution of two functions $d(f*g)=f(0)g+g*df$, where there is a jump part in $g$, and $dk(t)=-\omega k(t)dt$. We can use a stochastic differential equation(SDE) to describe $x(t)$ (c.f.\citep{Wang2016A}).

\begin{equation}\label{main}
dX_i=\omega(b_i-X_i(t))dt+ \sum_j \alpha h(X_j,X_i)dN_j
\end{equation}

It's  also necessary to add $\sigma dW$ into the model to represent noise, where $W$ is the Weiner noise with $W[0]=0$.
\begin{equation}
dX_i =\omega(X_i(0)-X_i(t))dt+\sum_{j} \alpha h(X_j,X_i)dN_j+ \sigma dW_i
\end{equation}

The above two equations are given in \citep{Wang2016A}. Because of the dependence on the path of the history, it is apparently not a Markov process.

We consider the limit that the population size tends to infinity. Let
 $X^N=(X_1,X_2\cdots, X_N)$ be the opinion variable of $N$-particles system. $X$ satisfies
\begin{equation}
\begin{split}
dX_i =\omega(b_i-X_i(t))dt+\frac{1}{N}\sum_{j\neq i}^N \alpha h(X_j,X_i)dN_j+ \sigma dW_i
\end{split}
\end{equation}
where $b_i=X_i[0]$.

In the case that $\omega=0$ and $h(x,y)=(x-y) 1_{[0,\epsilon]}(|x-y|)$, this is a Deffuant-Weisbunch like model.

In this paper, we will show that, when $N$ tends to infinity, there is a limit process $X$ for $X^N$, with SDE
\begin{equation}
\begin{aligned}
dX=\int dP(y) h(y,X)\lambda dt +\omega(b-X)dt +\sigma dW
\end{aligned}
\end{equation}
where $P$ is the law of $X$, $b=X[0]$ is a random variable.

So we can see that, given the initial distribution $\mu$, the law of $X$ is $P=\int P^b d\mu(b)$, where $P^b$ is the solution of the Fokker-Planck equation associated with the Mckean-Vlasov process:
\begin{equation}
\frac{\partial}{\partial t} P^b =-\frac{\partial}{\partial x} [\beta^{b,P} P^b]+\frac{1}{2} \sigma^2\frac{\partial^2}{\partial x^2}P^b
\end{equation}
$$P^b(0)=\delta_b$$
$$\beta^{b,P}=\int \mu(db')\int P^{b'} (dy) \alpha \lambda h(y,x)+\omega(b-x)$$
where $\delta_b$ is a Dirac measure with support on b.\\
Stemming from the physical meaning, we can assume that the support of initial distribution $\mu(b)$ is bounded, such as the uniform distribution on [0,1]. Also, for the proof of Mckean-Vlasov limit, we need the function $h(x)$ to have good enough properties to make
\begin{equation}\label{A}
E[h(X^{b'}_j,X^{b}_i)]< \infty
\end{equation}
and
\begin{equation}\label{B}
\int x^2 P(dx) <\infty
\end{equation}
where $P$ is the solution of the Mckean-Vlasov equation. If $h$ is bounded, the first condition is satisfied trivially.
We will give other examples of $h$ that can make these conditions be satisfied in the appendix. Also, we assume that $h(x,y)$ is Lipschitz for both of the variables(However, this is not the necessary condition c.f. \citep{M1996Asymptotic})
 \begin{equation}
 |h(x_1,y_1)-h(x_2,y_2)|<L|x_1-x_2|+K|y_1-y_2|
 \end{equation}
 Most of these conditions can be relaxed, but we hope to prove our results without excessive technical details on tightness, therefore  these restrictions are imposed.

\subsection{The outline of the proof}
\subsubsection{SDE related to initial values}
Since the process we study is not a Markov process, we need to consider SDE with the form
\begin{equation}
dX=F(X[0],X,t)dt+dW
\end{equation}
It has no  Markov generator. But the non-Markovian of it is "not so bad". The law of $X$ can be solved by the following method. For a given $X[0]$, we can solve the law by $e^{tA} P_0 $, where $A$ is the generator of the Fokker-Planck equation and $P_0=\delta_{x_0}$ is the initial Dirac distribution . Then the solution of the SDE has the form£¬
\begin{equation}\label{trick}
\int dP_0(X_0) e^{tA(X_0)} \delta_{X_0}
\end{equation}
where $A(X_0)$ is the generator with given $X[0]$, $\delta_{X_0}$ is the Dirac measure with the support on $X_0$ and $P_0$ is the initial distribution of $X_0$.

The essence of this method is to regard $X[0]$ as an additional random variables. When SDE has the above form, we need the double layer empirical measure, which will regain the symmetry.

We will consider the SDEs
\begin{equation}
dX_i =\omega(b_i-X_i(t))dt+\sum_{j\neq i} \alpha h(X_j,X_i)dN_j+ \sigma dW_i
\end{equation}
where $b_i$ is a stochastic variable such that $X_i[0]=b_i$, then we will prove the large number law for the double layer empirical measure$\frac{1}{N}\sum_i \delta_{X_i,b_i}$, which will give the Mckean-Vlasov limit.
This is equal to consider SDEs
\begin{equation}
dX^b_i =\omega(b-X^b_i(t))dt+\sum_{j\neq i} \int d\mu(b')\alpha h(X^{b'}_j,X^{b}_i)dN_j+ \sigma dW_i
\end{equation}
where $X^b$ is the stochastic process with given b and $E(X)=E(\int d\mu(b) X^b)$

\subsubsection{Intermediate process}

In \citep{Andreis2017McKean}, the authors proved that the Mckean-Vlasov limit of the following equation
\begin{equation}
dX_i ^N=F(X_i^N)dt +\sigma(X_i^N)dW_i+\frac{1}{N} \sum_{j\neq i}^N h(X_j^N,X_i^N)dN_j
\end{equation}
has the form
\begin{equation}
dX=F(X)+\int dP(Y) h(Y,X) \lambda dt+\sigma(X) dW
\end{equation}
where P is the law of X.

We use the same way as \citep{Andreis2017McKean} to reduce the jump-SDE to the "averaging dynamics". Their method is to decompose the jump terms into a martingale and a continuous part, then the martingale can be easily estimated by Burkholder-Davis-Gundy inequality.

We will show that, there is a intermediate process $Y_i$ with following SDE
\begin{equation}\label{inter}
dY^b_i=(b-Y_i)dt +\sigma dW_i +\int d\mu(b') \frac{1}{N}\sum_{j=1}^N \alpha h(Y^{b'}_j,Y^b_i)\lambda dt
\end{equation}
such that
\begin{equation}\label{result1}
\int d\mu(b)\frac{1}{N}\sum_{j=1}^N E[\sup||X^b_i-Y^b_i||]\leq \frac{C}{\sqrt{N}}.
\end{equation}

By comparing these two process directly, we can notice that
\begin{equation}
\int \mu(b) E[\sup_{r}||X_i^b-Y_i^b||]\leq \int \mu(b) (F+\Theta)
\end{equation}
where $F=\omega\int_0^t ||X_i^b(s)-Y_i^b(s)|| ds$\\
$\Theta=E[\sup_r{||\frac{1}{N}\sum_j\int_0^r ds (\alpha h(X_j,X_i)dN_j(s) -\alpha h(Y_j,Y_i)\lambda ds)||}]$\\
Using Doob Meyer decomposition
$$dN=d\widetilde{N} +\lambda dt$$
where $\widetilde{N}$ is a martingale.

Intuitively, the supremum of a martingale would not be very large. This can be proved by Burkholder-Davis-Gundy inequality and our assumption on $h$. Then we can obtain the approximation(\ref{result1}). See the second part of appendix for a detailed proof of (\ref{result1}).

\subsubsection{Mckean-Vlasov process}
The intermediate process has the form like
\begin{equation}
dX_i=\frac{1}{N}\sum_j f(X_i,X_j)dt +g(X_i,\omega_i)dt+dW
\end{equation}
where $\omega_i$ is a random variable.

The asymptotic behavior of the double layer empirical measure, with the case that $f$ and $g$ are bounded, has been studied in \citep{Pra1996McKean}. They use Varadhan's lemma to reduce the system to the case without interaction. Since f and g are not bound, we use the method in \citep{Sznitman1991Topics}to compare the two SDEs directly.
Using the Lipschitz condition, we can control $E[\sup_r||X(r)-\bar{X}(r)||]$ by the first Wasserstein distance $\rho$.

The first Wasserstein distance is defined by
$$W_1(\mu_1,\mu_2)=\inf{\int_{R^d\times R^d} |x-y| \eta(dx,dy): \eta\in H(\mu_1,\mu_2)}$$
where $H(\mu_1,\mu_2)$ is the set of all probability measures on $R^d \times R^d$ with marginals $\mu_1$ and $\mu_2$.

By Kantorovich-Rubinstein duality theorem, it has a dual representation:
$$\sup_{||f||_L\leq 1}{\int f(x) d(\mu_1-\mu_2)(x)}$$
where $||f||_L$ is the minimal Lipschitz constant for $f$.

Using this representation, since we have assumed that the coefficients are Lipschitz, we can see that
\begin{equation}
E[\sup_r||X(r)-\bar{X}(r)|| \leq L\int_0^T E[\rho(\mu^N_X(t),\mu_t)]dt
\end{equation}
where $\mu^N_X$ is the empirical measure, $\mu_t$ is the solution of Mckean-Vlasov equation.

For $E[\rho(\mu^N_X(t),\mu_t)]$, we have
\begin{equation}
E[\rho(\mu^N_X(t),\mu_t)]\leq E[\rho(\mu^N_X(t), \mu^N_{\bar{X}}(t))]]+E[\rho(\mu^N_{\bar{X}}(t),\mu^N_t]]
\end{equation}
where $\mu^N_{\bar{X}}$ is the empirical meausre for N independent Mckean-Vlasov processes.

Then we can derive our result by considering the Wasserstein distance between the solution of Mckean-valsov equation and empirical measure. Using  theorem 1 in \citep{Fournier2013On}, if the two order moment is finite:
$$E[\bar{X}^2]< \infty$$
when the size trends to infinity, the Wasserstein distance trends to zero, so our clam follows. The detailed proof is given in the appendix.

\section{The Characteristics of Mckean Vlasov dynamics}

\subsection{Mckean-Vlasov equation}
For the equation
\begin{equation}
\frac{\partial}{\partial t} P^b =-\frac{\partial}{\partial x} [\beta^{b,P} P^b]+\frac{1}{2} \sigma^2\frac{\partial^2}{\partial x^2}P^b
\end{equation}
$$P^b(0)=\delta_b$$
$$\beta^{b,P}=\int \mu(db')\int P^{b'} (dy) \lambda h(y,x)+\omega(b-x)$$
we consider a simple case

$$\lambda=1$$
$$\sigma=1$$
$$h(y,x)=\alpha (y-x)$$
Then we can see that
$$\beta^{b,P}= \alpha (m-x) +  \omega(b-x)  $$
where $m=\int xdP(x)$, the mean value of x.

Since h is an odd function, it is easy to see that
\begin{equation}
\begin{aligned}
\frac{dE(X(t))}{dt}=E( [\alpha \int P(dY)(Y-X) +\omega (b-X(t))])=\\
0+\omega E[X(0)-X(t)])
\end{aligned}
\end{equation}
So $E(X(t)-X(0))=Ce^{-\omega t}$, since $E(X(0)-X(0))=0$, $C=0$, $m=E(X(t))=E(X(0))$.

Then we have the Fokker-Planck equation associated with the Mckean-Vlasov process:
\begin{equation}
\frac{\partial}{\partial t} P^b =-\frac{\partial}{\partial x} [\alpha (m-x) +  \omega(b-x) P^b]+\frac{1}{2} \frac{\partial^2}{\partial x^2}P^b
\end{equation}

Integral over $\mu(b)$
\begin{equation}\label{me}
\frac{\partial}{\partial t} P =-\frac{\partial}{\partial x} [(\alpha(m-x))P + \omega \int d\mu(b) (b-x) P^b]  +\frac{1}{2}\frac{\partial^2}{\partial x^2}P
\end{equation}

Another example  is the  bounded confidence type model as in \citep{Hegselmann2005Opinion,Guillaume2000Mixing}, where every agent interacts
only within a certain level of confidence.
\begin{equation}
h(y,x)=(y-x)k(y-x)
\end{equation}
where $k(y-x)$ is a continuous function and if $y-x <\Delta_1$,$k(y-x)=1$, and if $y-x >\Delta_2 $,$k(y-x)=0$. Since in this case, $h(x)$ is bounded and Lipschitz, our assumptions are satisfied.

Note that for the case that $h(y,x)=\alpha y$ which is used in SLANT system in\citep{De2015Learning},
$$\frac{\partial}{\partial t}P^b=-\frac{\partial}{\partial x} [\alpha m(t) +  \omega(b-x) P^b]+ \frac{1}{2}\sigma^2 \frac{\partial^2}{\partial x^2}P^b$$
where $m(t)=\int x dP(x,t) $. Since
\begin{equation}
\begin{aligned}
\frac{dE(X(t))}{dt}&=E( [\alpha \int P(dY) Y +\omega (b-X(t))])\\
&=(\alpha-\omega)E(X)+\omega E(X(0))
\end{aligned}
\end{equation}
As a sanity check, this formula matches equation(4) in \citep{De2015Learning}.
Solving this equation, we get:
\begin{equation}\label{m}
\begin{aligned}
m(t)+\frac{\omega}{\alpha-\omega}m(0)=C e^{(\alpha-\omega)t}\\
\end{aligned}
\end{equation}
\begin{equation}
C=\frac{\alpha}{\alpha-\omega}
\end{equation}

\subsection{ Comparison with other models }
Ref.\citep{Como2011SCALING} showed that when the population size tends to infinity, the limit behavior of Deffuant-Weisbunch model can be described by a measure-valued ODE. And also, in \citep{Toscani2006Kinetic} Toscani introduced a kinetic model of opinion formation. In this model, the opinion is changed by the binary interaction of collision, and, the dynamics of the opinion distribution is modelled by Boltzmann type integro-differential equation.

In these models, the dynamics of opinion distribution is depicted as follows:
\begin{equation} \label{1}
\frac{d}{dt}<\phi,\mu_t> =<\phi,H(\mu_t)>
\end{equation}
where $\mu_t$ is the measure of opinion distribution and $\phi$ is an arbitrary test function.
The right side of (\ref{1}) has the following form \citep{Como2011SCALING}
$$<\phi,H(\mu_t)>=\int \int (\phi(1-\omega )x+\omega y)-\phi(x))k(x,y) d\mu(x) d\mu(y)$$

And in Boltzmann case \citep{Toscani2006Kinetic}
\begin{equation}
 <\phi,H(\mu_t)>=\int \int \phi(w')+\phi(w'_*)-\phi(w)-\phi(w_*) d\mu(w) d\mu(w_*)
\end{equation}
$$w'=(1-\eta P(w,w_*))w+\eta P(w,w_*)w_*$$
$$w'_*=(1-\eta P(w_*,w))w+\eta P(w_*,w)w$$

Intuitively, both models can be derived from the stochastic differential equation with Poisson jump interactions by the method mentioned in \citep{M1996Asymptotic}. Events that interact with other particles are regarded as a Poisson process.

For the  system with $h(y,x)=y-x$:
$$\frac{\partial}{\partial t} P =-\frac{\partial}{\partial x} [(\alpha(m-x))P + \omega \int d\mu(b) (b-x) P^b]  +\frac{1}{2}\frac{\partial^2}{\partial x^2}P$$
in the case that $\omega=0$, which is the non-self-excitation case,it  degenerates into the general opinion  dynamics model, or equally, the quasi-invariant opinion limit for Toscani's Boltzmann equation in \citep{Toscani2006Kinetic}:
\begin{equation}
\frac{\partial}{\partial t}P=\frac{\lambda}{2}\frac{\partial^2}{\partial x^2}(D(|x|^2)P) -\frac{\partial}{\partial x}(P(|x|)(m-x)P)
\end{equation}

The $\omega$ term is   a modification to consider the influence of self excitation process.

\subsection{Steady state distribution}
Firstly, we consider a simple case that the initial distribution $P_0=\frac{1}{2}\delta_{-10}+\frac{1}{2}\delta_{10}$. We will use this simple form to explain how the initial distribution is coupled to the equation. This distribution is to represent two groups of people. Half of them lay on one side and the others is on the contrary. Then the equation will become
\begin{equation}
\begin{aligned}
\frac{\partial}{\partial t}P^{\{-10\}}= -\frac{\partial }{\partial x} [(\alpha(\bar{x}-x))P^{\{-10\}} +\\
 \omega (-10 -x) P^{\{-10\}}]  +\frac{1}{2} \frac{\partial^2}{\partial x^2}P^{\{-10\}}
\end{aligned}
\end{equation}
And
\begin{equation}
\frac{\partial}{\partial t} P^{\{10\} } =-\frac{\partial}{\partial x} [(\alpha(\bar{x}-x))P^{\{10\}} + \omega(10-x) P^{\{10\}}]  +\frac{1}{2}\frac{\partial^2}{\partial x^2}P^{\{10\}}
\end{equation}
The  overall distribution $P$ is equal to $\frac{1}{2}P^{\{-10\}}+\frac{1}{2}P^{\{10\}}$.
\begin{figure*}
  \includegraphics[width=0.45\textwidth]{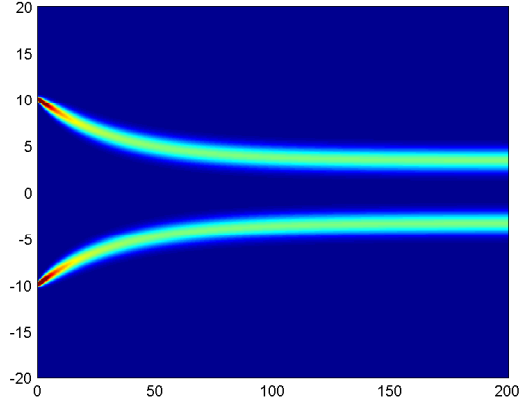}
\caption{Evolution of the opinion density,with $\sigma =0.02 ,\alpha=0.02,\omega=0.01$.}
\label{fig:1}       
\end{figure*}

\begin{figure*}
\begin{minipage}[t]{0.5\linewidth}
\centering
\includegraphics[width=3.5in]{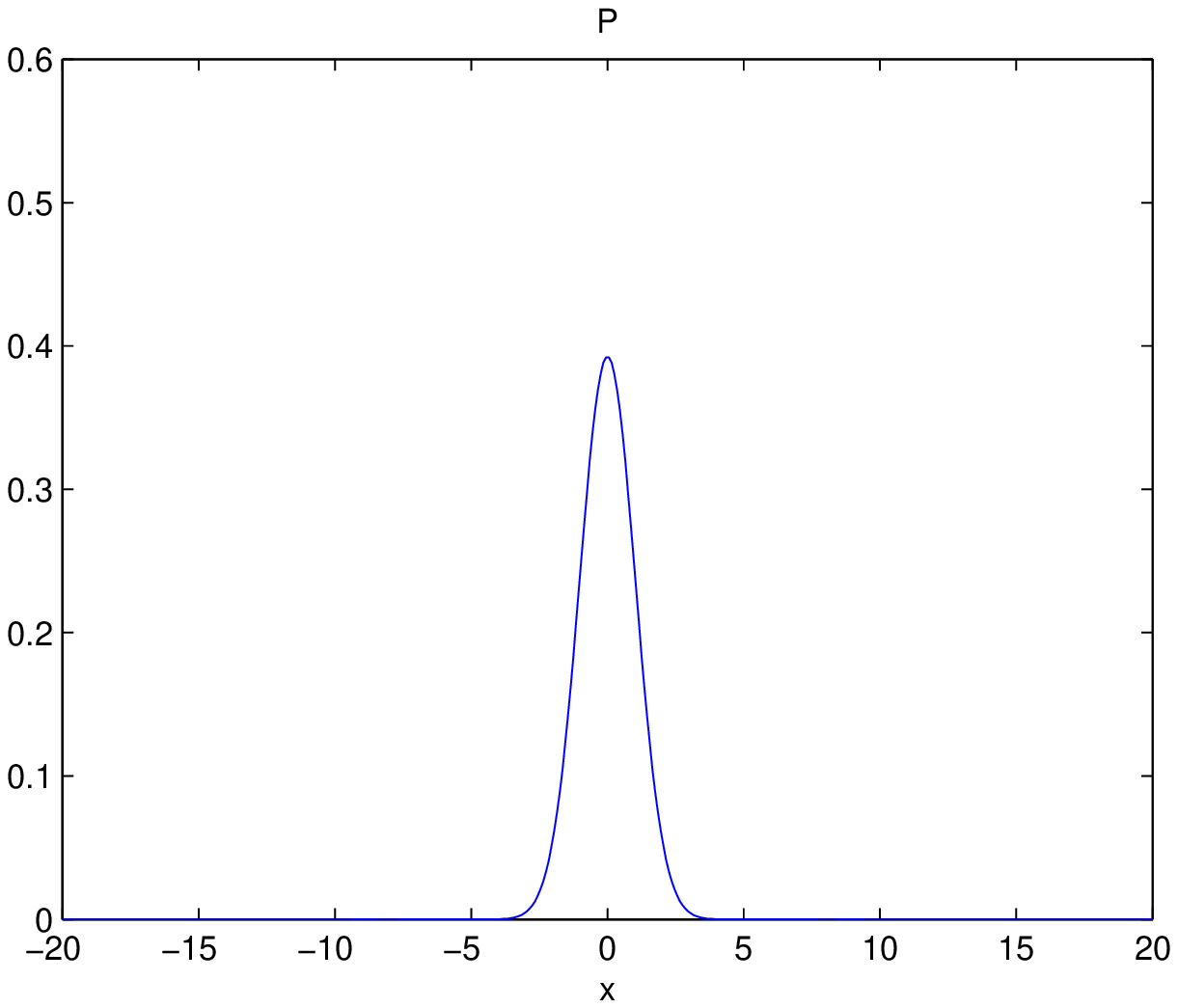}
\caption{Steady state for the dynamics with parameters\\ $\sigma=0.02, \alpha=0.02,\omega$=0}
\label{fig:2}
\end{minipage}%
\begin{minipage}[t]{0.5\linewidth}
\centering
\includegraphics[width=3.5in]{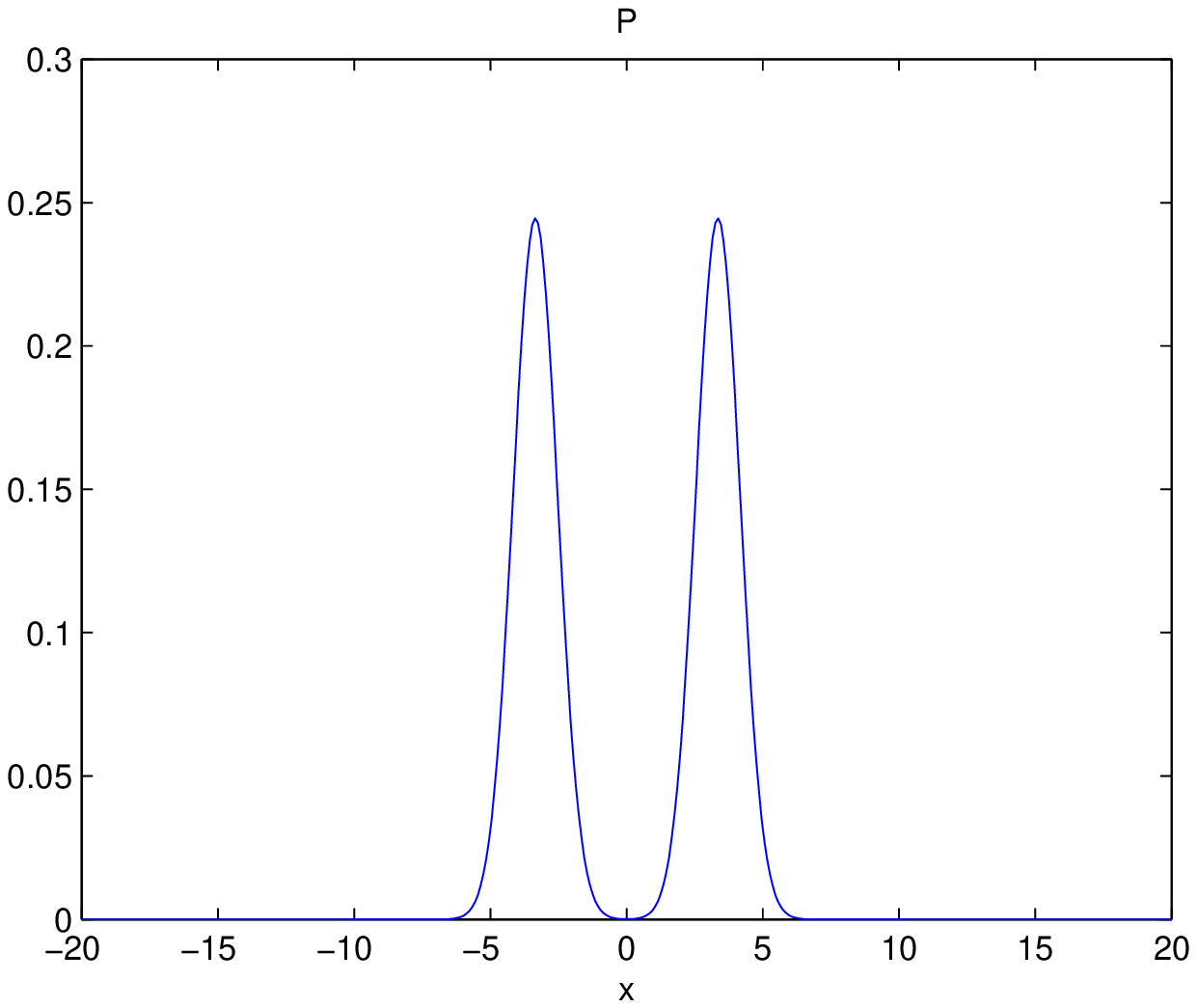}
\caption{Solution profiles at time t=200, with parameters\\ $\sigma=0.02, \alpha=0.02,\omega$=0.01}
\label{fig:3}
\end{minipage}
\end{figure*}

The evolution process of this distribution is shown in Figure \ref{fig:1}. And Figure \ref{fig:2} and Figure \ref{fig:3} show the differences of two kinds of steady states.

In the case that $\omega=0$, which is a Markov process, the Mckean-Vlasov equation has the form:
$$\frac{\partial}{\partial t} P =-\frac{\partial}{\partial x} (\alpha(m-x))P +\sigma^2\frac{1}{2}\frac{\partial^2}{\partial x^2}P$$

Such an equation has a steady state
$$-\frac{\partial}{\partial x} (\alpha(m-x))P +\sigma^2\frac{1}{2}\frac{\partial^2}{\partial x^2}P=0$$
And this steady state is asymptotically stable, such that when $t\to \infty$, the solution of the equation will trend to the steady state. It is easy to see that the steady state is a normal distribution as Figure \ref{fig:2}.

However, this $\omega=0$ model can not describe the actual situation, since in this model, the final state of the opinion distribution have no community structures but a perfect consensus, which is so called Abelson¡¯s diversity puzzle \citep{Abelson1964Mathematical}, a persistent research puzzle in the social sciences. Generally, this puzzle is solved by the bounded-confidence mode, such as HK\citep{Hegselmann2005Opinion} and DW\citep{Guillaume2000Mixing} models. It is proved in \citep{Como2011SCALING} for the large number of individuals case and \citep {Zhang2012Convergence} for the general case, that the steady state of DW model has the form $\sum_i C_i\delta_{x_i}$, where $x_i$ is the opinion of the community i, and $|x_i-x_j|>R$ if $i\neq j$, where R is the bounded confidence distance.

The equation with $\omega>0$ has the steady state as Figure \ref{fig:3} when the initial distribution has the form $\frac{1}{2}(\delta_{-10}+\delta_{10})$. In the general case, let $P=\int d\mu(b) P^b$, and $P^b$ is the solution of Mckean Vlasov equation
\begin{equation}
\frac{\partial}{\partial t} P^b =-\frac{\partial}{\partial x} [\alpha (m-x) +  \omega(b-x) P^b]+ \frac{1}{2}\sigma^2 \frac{\partial^2}{\partial x^2}P^b
\end{equation}
The steady state of $P^b$ is  easy to calculate:
\begin{equation}\label{st}
P^b(\infty)=\frac{1}{Z}e^{-(\frac{(\alpha +\omega)x-(\alpha m+\omega b)}{\sigma})^2}
\end{equation}
where Z is the  Normalization constant.

So we have $P(\infty)=\int d\mu(b) P^b(\infty)$.

Considering the limit behavior $\sigma\to 0$, the steady state of Mckean-Vlasov equation when $h(y,x)=y-x$ has the form $\int d\mu(b) \delta(x-\frac{\alpha m +\omega b}{\alpha +\omega})$. So that even in the case that the variance of the noise trends to 0, this equation will not converge to a single Dirac measure, which is different from the $\omega=0$ case. Intuitively, the term $\omega(b-x)$ make the opinion try not to deviate from the initial position too far, which can be considered as a soft bounded confidence.

In the case that $h(y,x)=\alpha y$ which is used in SLANT system in\citep{De2015Learning}, supposing $\omega>\alpha$, the steady state will have the similar form since the equation of which is
$$-\frac{\partial}{\partial x} [\alpha m +  \omega(b-x) P^b]+ \frac{1}{2}\sigma^2 \frac{\partial^2}{\partial x^2}P^b=0$$
where $m=m(\infty)=\frac{\alpha}{\omega-\alpha}m(0) $(see eq \ref{m} ), the final mean value of x. We can show that  $P(\infty)=\int d\mu(b) \delta(x-\frac{\alpha m +\omega b}{\omega})$ easily.

\section{Conclusion}
In this paper, we study the asymptotic behaviors of the  exponential function weight self excitation interaction system $X^N$ on the large homogeneous network when $N$ the number of the individuals trends to infinity. We prove that there is a Mckean-Vlasov  process $\bar{X}$, such that $\lim \limits_{N \to +\infty} E|X^N-\bar{X}|$ is zero, whose opinion distribution evolve according to a Mckean-Vlasov type integro differential equation which couples with the initial distribution. The steady state of this PDE is also studied. We show that even if we do not consider bounded confidence, it will give not a perfect consensus, but a distribution  coupling with the initial distribution, such that this model can avoid  Abelson¡¯s diversity puzzle.

\section*{Acknowledgements}
The authors thank Wenyu Zhang for his great help in writing and we thank Prof. Haibo Wang who introduced the recent developments in the mean field game theory to us.\\
This research is supported by the National Key Research and Development Program of China (No.2018YFC0831306), the Fundamental Research Funds for the Central Universities (No.2017JBZ107, No.2018YJS003) and the National Natural Science Foundation of China under grant 61271308.

\appendix

\begin{appendix}

\section{The condition under which the assumptions hold}
We want to give some examples that the assumptions can be satisfied. We will show that when $h$ has linear growth , the assumptions (\ref{A},\ref{B}) are satisfied.

For  SDE
\begin{equation}
dX^b_i=\frac{1}{N}\sum_{j\neq i} \int d\mu(b') \alpha h(X^{b'}_j,X^b_i)dN_j +\omega(b-X^b_i)dt+\sigma dW_i
\end{equation}
it is shown in \citep{Situ2005Theory} that since both $b-X$ and $h$ has linear growth, this SDE has a strong solution, and also we can estimate the moment.

In order to prove (\ref{A}), since h is linear growth, we can turn to consider the second moment by Jensen's inequality $E[h(y,x)]\leq L(E[|y|]+E[|x|])\leq 2L(E[x^2])^{\frac{1}{2}}$. Using Ito formula,
\begin{equation}
\begin{aligned}
\int d\mu(b) \frac{1}{N}\sum_i (X^b_i(T))^2 =\\
\int d\mu(b) [\frac{1}{N}\sum_i (X^b_i(0)^2 + \int_0^T X^b_i(s) b -(X^b_i(s))^2 ds )]\\
+\int d\mu(b) d\mu(b')\frac{1}{N} \sum_{i,j} [ \int_0^T 2X_i^b(s) \frac{1}{N}\alpha h(X_j^{b'},X_i^b)\\
+ \frac{1}{N^2}\alpha^2 |h(X_j^{b'},X_i^b)|^2 dN_j ]\\
+ C+martingale
\end{aligned}
\end{equation}

Since $|h(x,y)|\leq L(|X|+|Y|)$ and $dN=\lambda dt + martingale$ we can show that
\begin{equation}
\begin{aligned}
\int d\mu(b) \frac{1}{N}\sum_i E(X^b_i(T))^2 &\leq \int d\mu(b) \frac{1}{N}\sum_i[(X^b_i(0)^2 \\
&+K\int_0^T E(X^b_i(s))^2+C ds \\
&+\int_0^T E(X^b_i(s)) b -E(X^b_i(s))^2 ds )]
\end{aligned}
\end{equation}

We use $E[X]\leq E[|X|]\leq (E[X^2])^{\frac{1}{2}}\leq K(1+E[X^2])$ by Jensen's inequality. Then using Gronwall' lemma, our clam follows.

As for the second conditon(\ref{B}), the above method can be also applicable.

\section{Reduce to the averaging dynamics}

The SDE for stochastic process $X_i$ is
\begin{equation}
dX^b_i=g(b,X^b_i)dt+\sigma dW_i+ \int d\mu(b') \frac{1}{N}\sum_{j\neq i}^N \alpha h(X^{b'}_j,X^b_i)dN_j
\end{equation}

Following \citep{Andreis2017McKean}, it is useful to introduce a intermediate process $Y_i$ with SDE
\begin{equation}
dY^b_i=g(b,Y^b_i)dt+\sigma dW_i+  \int d\mu(b')  \frac{1}{N} \sum_{j=1}^N \alpha h(Y^{b'}_j,Y^b_i)\lambda dt
\end{equation}
Since the Possion process is a semimartingale, $N(t)$ can  be decomposed into $N(t)=\widetilde{N}(t)+\lambda t$ by Doob-Meyer decomposition, where $\widetilde{N}(t)$ is a martingale. The $\lambda dt$ in Y is the second part of Doob-Meyer decomposition of N. To simplify the symbol, we set $\lambda=1$

We will show that, for large N,
\begin{equation}
\int d\mu(b) \frac{1}{N}\sum_{j=1}^{N} E[\sup\limits_{r\in [0,t]}||X^b_j(r)-Y^b_j(r)||]\leq \frac{C}{\sqrt{N}}
\end{equation}

so X can be approximated by Y.

Let
\begin{equation}
G^b_i=E[\int_0^t||g(X^b_i,b)-g(Y^b_i,b)||ds]
\end{equation}
\begin{equation}
\begin{aligned}
\Theta_i=E[\sup_r||d\mu(b')[\frac{1}{N}\int \sum_{j\neq i} \int_0^r \alpha h(X^{b'}_j,X^b_i)dN_j-\\
\frac{1}{N}\sum_j \int_0^r \alpha h(Y^{b'}_j,Y^b_i) ds]||]
\end{aligned}
\end{equation}

So
\begin{equation}
\begin{aligned}
\int d\mu(b) \frac{1}{N}\sum_i E[\sup_r ||X^b_i(r)-Y^b_i(r)||]\leq \\
\int d\mu(b) [\frac{1}{N}\sum_i (G_i+ \Theta_i)]
\end{aligned}
\end{equation}
In our model (\ref{main}), $g(x_i,b)=\omega (b-x_i)$, so that we have $G_i=E[\int_0^t ||\omega(Y^b_i-X^b_i)||]$.

For $\Theta$
\begin{equation}
\begin{aligned}
\Theta_i &\leq E[\sup_r||\int d\mu(b')  \frac{1}{N}\sum_{j \neq i} \int_0^r \alpha h(X^{b'}_j,X^b_i) d \widetilde{N}(t)||]\\
&+E[\sup_r||\int d\mu(b') \frac{1}{N}\sum_j \int_0^r \alpha(h(X^{b'}_j,X^b_i)-h(Y^{b'}_j,Y^b_i)) ds||]\\
&+\frac{1}{N}E[\sup_r ||\int_0^r \int d\mu(b')\alpha h(X^{b'}_i,X^{b}_i) ds||]
\end{aligned}
\end{equation}
As in \citep{Andreis2017McKean}, we can use the Burkholder-Davis-Gundy inequality for martingales.
\begin{equation}
\begin{aligned}
\Theta_i &\leq \frac{k}{N}  E[( \sum_{j\neq i} \int_0^t||\int d\mu(b') h(X^{b'}_j,X^b_i)||^2 ds)^{1/2}]\\
&+ \frac{1}{N}\int_0^t E[\sup_r  \sum_j \alpha || \int d\mu(b') h(X^{b'}_j,X^b_i)-h(Y^{b'}_j,Y^b_i)|| ]ds +\frac{H}{N}\\
&\leq \frac{k}{N} \int_0^t E[(  N ((|| \int d\mu(b') h(X^{b'}_j,X^b_i)|| ds)^2 )^{1/2}]\\
&+ \frac{1}{N} \int_0^t  E[  \sum_j \alpha ||\int_0^r \int d\mu(b') h(X^{b'}_j,X^b_i)-h(Y^{b'}_j,Y^b_i)|| ]ds +\frac{H}{N}\\
&\leq \frac{k}{N}  E[(  N t((\sup_r{||\int d\mu(b') h(X^{b'}_j(r),X^b_i(r))||} )^2 )^{\frac{1}{2}} ]\\
&+ \int_0^t  \frac{1}{N}E[ \sum_j \alpha ||\int d\mu(b') h(X^{b'}_j,X^b_i)-h(Y^{b'}_j,Y^b_i)|| ]ds  +\frac{H}{N} \\
&\leq \frac{C}{\sqrt{N}}+L\int_0^t E[||X^{b}_i-Y^b_i||]ds +\\
&\frac{L}{N}\sum_j \int d\mu(b') \int_0^t E[||X^{b'}_j-Y^{b'}_j||]ds +\frac{H}{N}
\end{aligned}
\end{equation}
where we use the symmetry, assumption (\ref{A}),and the Lipschitz condition for h.
Combine the two results, then we get
\begin{equation}
\begin{aligned}
&\int d\mu(b) \frac{1}{N}\sum_j E[\sup{||X^b_i-Y^b_i||}\leq\\
&\frac{C}{\sqrt{N}}+K \int d\mu(b) \frac{1}{N}\sum_j \int_0^t E[||X^b_i-Y^b_i||]ds +\frac{H}{N}\\
&\leq \frac{C}{\sqrt{N}} +K \int d\mu(b)  \frac{1}{N}\sum_j \int_0^t E[\sup{||X^b_i-Y^b_i||}]ds +\frac{H}{N}
\end{aligned}
\end{equation}
Then using the standard technology of Gronwall's lemma, we get the approximation.

\section{Mckean-Vlasov process}
As shown in the previous section,we can turn to analyze the stochastic process Y.
Let $\bar{X}$ be a stochastic process with SDE:
\begin{equation}\label{MC}
d\bar{X}=\int dP(y) \alpha h(y,\bar{X}) dt +\omega(b-\bar{X})dt +\sigma dW
\end{equation}
where $P$ is the law of $\bar{X}$, and b is a  random variable such that $b=X[0]$.\\
 Using the trick in (\ref{trick}), we can show that
\begin{equation}
P(y)=\int d\mu(b) P^b(y)
\end{equation}
where $P^b$ is the law of SDE with given b, the initial distribution is $\delta_b$, and $\mu(b)$ is the given initial distribution of $\bar{X}$.
We will show that
\begin{equation}
\lim_{N\to \infty} \int d\mu(b) \frac{1}{N}\sum_i E[\sup_{t\leq T}|Y_t^{b,i}-\bar{X}_t^{b,i}|] =0
\end{equation}
Let $Q^N$ be the empirical measures $\sum_j\frac{1}{N} \delta_{\bar{X_j}}$, we have
\begin{equation}
\begin{aligned}
Y_t^{b,i}-\bar{X}_t^{b,i}&=\\
&\int_0^t \frac{1}{N} \int d\mu(b') \sum_j h(Y_s^{b,i},Y_s^{b',j}) \\
&- \int h(\bar{X}_s^{b,i},y) Q(dy)ds \\
&+\omega \int_0^t \bar{X}_s^{b,i}-Y_s^{b,i} ds \\
\end{aligned}
\end{equation}
\begin{equation}
\begin{aligned}
&=\int_0^t ds \frac{1}{N} \int d\mu(b') \sum_j h(Y_s^{b,i},Y_s^{b',j})\\
&-h(\bar{X}^{b,i}_s,Y^{b',j}_s)\\
&+(h(\bar{X}^{b,i}_s,Y_s^{b',j})-h(\bar{X}_s^{b,i},\bar{X}_s^{b',j}))\\
&+\frac{1}{N}\sum_j \int d\mu(b') (h(\bar{X}^{b,i}_s,\bar{X}_s^{b',j})\\
&-\int h(\bar{X}_s^{b,i},y)Q(dy))\\
&+\omega \int_0^t \bar{X}_s^{b,i}-Y_s^{b,i} ds\\
\end{aligned}
\end{equation}

Since h is Lipschitz,  we can see that
\begin{equation}
\begin{aligned}
\int d\mu(b) &E[|\sup_r{Y^{b,i}(r)-\bar{X}^{b,i}(r)|}]\leq\\
&L\int d\mu(b) \sup_r \int_0^r ds (E[|Y^{b,i}-\bar{X}^{b,i}|]\\
&+\frac{1}{N}\sum_j \int d\mu(b') E[|Y^{b',j}-\bar{X}^{b',j}|]\\
&+\int d\mu(b) E|\frac{1}{N}\sum_j \int d\mu(b') h(\bar{X}^{b,i},\bar{X}^{b',j})\\
&-\int h(\bar{X}^{b,i},y)Q(dy)]\\
&+\omega \int_0^t ||\bar{X}_s^{b,i}-Y_s^{b,i}||\\
\end{aligned}
\end{equation}\\
Summing over i
\begin{equation}
\begin{aligned}
\frac{1}{N} \sum_i &\int d\mu(b) E[\sup_r {|Y^{b,i}(r)-\bar{X}^{b,i}(r)|}]\leq\\
&\frac{1}{N}\{ L\sup_r \int_0^r ds  \sum_i\int d\mu(b) [E[|Y^{b,i}-\bar{X}^{b,i}|]\\
&+E|\frac{1}{N}\sum_j\int d\mu(b') h(\bar{X}^{b,i},\bar{X}^{b',j})\\
&-\int h(\bar{X}^{b,i},y)Q(dy) \\
&+\omega \int_0^t \sum_i E||\bar{X}_s^{b,i}-Y_s^{b,i}||]\}
\end{aligned}
\end{equation}
then using Gronwall's lemma
\begin{equation}
\begin{aligned}
\frac{1}{N} &\sum _i \int d\mu(b) E[\sup_r {|Y^{b,i}(r)-\bar{X}^{b,i}(r)|}]\leq\\
&L(T)\sup_r \int_0^r \int d\mu(b)E|\frac{1}{N}\sum_j\int d\mu(b') h(\bar{X}^{b,i},\bar{X}^{b',j})\\&-\int h(\bar{X}^{b,i},y)Q(dy)] ds
\end{aligned}
\end{equation}
So we need to estimate the right side.\\

Using the Lipschitz condition for h again, this can be controled by $E[\rho (Q^N_{\bar{X}}(t),Q_{\bar{X}}(t))]$,where $\rho$ is the Wasserstein distance between the two measures, and $Q_{\bar{X}}$ is the solution of Mckean-Vlasov equation (\ref{MC}). $Q^N$ is the empirical measures $\sum_j\frac{1}{N} \delta_{\bar{X_j}}$.
It follows from  \citep{Fournier2013On} theorem 1 that if the moments is bound(This is satisfied by (\ref{B}) )
\begin{equation}
E[\bar{X}^2]< \infty
\end{equation}
our clam follows.\\
\end{appendix}
\\
 \bibliography{References}
 \end{document}